\documentclass{elsart}
\usepackage{amssymb,amsmath}
\usepackage{latexsym}

\begin{document}

\begin{frontmatter}

\title{On the Group Classification of Systems of Two Linear Second-Order
Ordinary Differential Equations with Constant Coefficients}

\author[SUT]{S.V. Meleshko},
\ead{sergey@math.sut.ac.th}
\author[DUT]{S. Moyo}
\ead{moyos@dut.ac.za}
\author[SUT]{G.F. Oguis}
\ead{fae@math.sut.ac.th}

\address[SUT]{Suranaree University of Technology, School of Mathematics,
Nakhon Ratchasima 30000, Thailand}
\address[DUT]{Durban University of Technology, Department of Mathematics,
Statistics and Physics \& Institute for Systems Science, P O Box 1334, Steve Biko Campus, Durban 4000, South Africa}
\address[DUT]{Suranaree University of Technology, School of Mathematics,
Nakhon Ratchasima 30000, Thailand}

\begin{abstract}
The completeness of the group classification of systems of two linear second-order ordinary differential equations with constant coefficients is delineated in the paper. The new cases extend what has been done in the literature. These cases correspond to the type of equations where the commutative property of the coefficient matrices with respect to the dependent variables and the first-order derivatives in the considered system does not hold. A discussion of the results as well as a note on the extension to linear systems of second-order ordinary differential equations with more than two equations are given.

\end{abstract}

\end{frontmatter}

\section{Introduction}
In this paper we consider the complete group classification of systems
of two linear second-order ordinary differential equations with constant
coefficients. Systems of second-order ordinary differential equations
arise in the modeling of physical phenomena in the areas of physics,
chemistry and mathematics. Of interest to the study of such systems
is their symmetry properties. The founder of symmetry analysis of
scalar ordinary differential equations, Sophus Lie \cite{bk:Lie[1883]},
gave a complete group classification of a scalar ordinary differential
equation of the form
\[
y^{\prime\prime}=f(x,y).
\]
 L.V.Ovsiannikov \cite{bk:Ovsiannikov[2004]}
 later performed Lie's
classification in a different way. This classification was obtained
by directly solving the determining equations and using the equivalence
transformations. The same approach was applied in \cite{bk:PhaukMeleshko[2013]}
for the group classification of more general types of equations of
the form $y^{\prime\prime}=P_{3}(x,y;y^{\prime})$, where $P_{3}(x,y;y^{\prime})$
is a polynomial of third degree with respect to the first-order derivative
$y^{\prime}$. In the general case of a scalar ordinary differential
equation $y^{\prime\prime}=f(x,y,y^{\prime})$ the application of
the method that involves directly solving the determining equations
gives rise to some difficulties. The group classification of such
equations \cite{bk:MahomedLeach[1989]} is based on the enumeration
of all possible Lie algebras of operators acting on the plane $(x,y)$.
Lie \cite{bk:Lie[1891b]} gave the classification of all dissimilar
Lie algebras (under complex change of variables) in two complex variables.
Later on, in 1992, Gonzalez-Lopez et al. ordered the Lie classification
of realizations of complex Lie algebras and extended it to the real
case \cite{bk:GonzalezKamranOlver[1992a]}. A large amount of results
on the dimension and structure of symmetry algebras of linearizable
ordinary differential equations is well-known (see
\cite{bk:MahomedLeach[1989],bk:MahomedLeach[1990],bk:HandbookLie_v3,bk:Gorringeleach[1988],%
bk:Ovsiannikov[1978],bk:WafoMah[2000],bk:BoykoPopovychShapoval[2012]}).

There are several papers where the group classification of systems
of second-order ordinary differential equations was considered. The
literature has dealt extensively with symmetry properties of a scalar
ordinary differential equation but the group classification of a system
of even two linear second-order equations with constant coefficients
is not complete. In recent works \cite{bk:BoykoPopovychShapoval[2012],bk:WafoSoh[2010],bk:Meleshko[2011],bk:Campoamor-Stursberg[2011],bk:Campoamor-Stursberg[2012]}
the authors focused on the study of systems of second-order ordinary
differential equations with constant coefficients of the form
\begin{equation}
{\bf y}^{\prime\prime}=M{\bf y},\label{feb27.1}
\end{equation}
 where $M$ is a matrix with constant entries. In the present paper
it is shown that these types of systems do not exhaust a set of all
systems of linear second-order ordinary differential equations with
constant coefficients.

In the general case of systems of two linear second-order ordinary
differential equations, the more advanced results are obtained in
\cite{bk:WafoMahomed[2000]}, where using the canonical form,
\[
{\bf y}^{\prime\prime}=\left(\begin{array}{cc}
a(x) & b(x)\\
c(x) & -a(x)
\end{array}\right){\bf y},
\]
 presented several representatives of nonequivalent classes. However
the list of all distinguished representatives of systems of two linear
second-order ordinary differential equations was not obtained in \cite{bk:WafoMahomed[2000]}.
Recently the complete group classification of linear systems of two
second-order ordinary differential equations has been studied in \cite{bk:MoyoMeleshkoOguis[2013]}.

Here the complete classification of linear systems of two second-order ordinary differential
equations with constant coefficients is presented.

As far as we are aware the results found here are new and have not
been reported in the literature.

The paper is organized as follows: In the first part we present the
simplification of a linear system of second-order equations with constant coefficients.
This is then followed by a preliminary study of linear systems of second-order
ordinary differential equations with constant coefficients. The later
part is followed by a complete treatment of linear systems of second-order
ordinary differential equations with constant coefficients as well as the discussion of the results and conclusion in the final part of the paper.

\section{Simplification of a system of linear equations}

Let a linear system of second-order equations with constant coefficients be given
by:
\begin{equation}
{\bf y}^{\prime\prime}=A{\bf y}^{\prime}+B{\bf y}+{\bf f},\label{eq:oct30.1}
\end{equation}
 where $A$ and $B$ are constant matrices. Using a particular solution
${\bf y}{}_{p}$, one can reduce (\ref{eq:oct30.1}) to the homogeneous
system
\begin{equation}
{\bf y}^{\prime\prime}=A{\bf y}^{\prime}+B{\bf y}.\label{eq:oct30.2}
\end{equation}

\subsection{Linear change of the dependent variables}

 Applying the change,
\[
{\bf y}=C{\bf x},
\]
 where $C=C(t)$ is a nonsingular matrix, system (\ref{eq:oct30.2})
becomes 
\begin{equation}
{\bf x}^{\prime\prime}=\bar{A}{\bf x}^{\prime}+\bar{B}{\bf x},\label{eq:oct30.3}
\end{equation}
 where
\[
\bar{A}=C^{-1}(AC-2C^{\prime}),\quad\bar{B}=C^{-1}(BC+AC^{\prime}-C^{\prime\prime}).
\]
 If one chooses the matrix $C(t)$ such that
\[
C^{\prime}=\frac{1}{2}AC,
\]
 then
\[
\bar{B}=C^{-1}\left(BC+\frac{1}{2}A^{2}C-\frac{1}{4}A^{2}C\right)=C^{-1}\left(B+\frac{1}{4}A^{2}\right)C.
\]
 The existence of the nonsingular matrix $C(t)$ is guaranteed by
the existence of the solution of the Cauchy problem
\[
C^{\prime}=\frac{1}{2}AC,\quad C(0)=E,
\]
 where $E$ is the unit matrix. This solution defines the matrix
\[
C(t)=e^{\frac{t}{2}A},
\]
 which commutes with the matrix $A$:
\[
AC=CA.
\]

Let us study the case when the matrix $\bar{B}$ is constant. First
of all notice that
\[
\frac{d}{dt}\left(C^{-1}\right)=-C^{-1}C^{\prime}C^{-1}.
\]
 Then one has
\[
\begin{array}{c}
\frac{d}{dt}\bar{B}=\frac{d}{dt}\left(C^{-1}\left(B+\frac{1}{4}A^{2}\right)C\right)\\
=-C^{-1}ACC^{-1}\left(B+\frac{1}{4}A^{2}\right)C+C^{-1}\left(B+\frac{1}{4}A^{2}\right)AC\\
=-C^{-1}A\left(B+\frac{1}{4}A^{2}\right)C+C^{-1}\left(B+\frac{1}{4}A^{2}\right)AC=C^{-1}(BA-AB)C.
\end{array}
\]
 Hence, the matrix $\bar{B}$ is constant if and only if the matrices
$A$ and $B$ commute
\begin{equation}
AB=BA.\label{eq:oct30.4}
\end{equation}

The papers cited in the literature \cite{bk:BoykoPopovychShapoval[2012],bk:WafoSoh[2010],bk:Meleshko[2011],bk:Campoamor-Stursberg[2011],%
bk:Campoamor-Stursberg[2012]},
where the group classification
of systems of linear equations with constant coefficients was considered,
studied systems with constant matrix $\bar{B}$. Thus, even for
systems of linear equations with constant coefficients there is no
complete study of the group classification. The present paper fills this niche.


\subsection{Example}

Let us choose the matrices
\[
A=\left(\begin{array}{cc}
0 & 1\\
1 & 0
\end{array}\right),\quad B=\left(\begin{array}{cc}
0 & 1\\
0 & 0
\end{array}\right),
\]
 then
\[
AB-BA=\left(\begin{array}{rc}
-1 & 0\\
0 & 1
\end{array}\right).
\]
 Notice that
\[
A^{2}=E.
\]
 Then
\[
\begin{array}{c}
e^{sA}=E+sA+\frac{s^{2}}{2!}A^{2}+\frac{s^{3}}{3!}A^{3}+\frac{s^{4}}{4!}A^{4}+...\\
=\left(\frac{s^{2}}{2!}+\frac{s^{4}}{4!}+...+\frac{s^{2n+1}}{(2n+1)!}...\right)E\\
+\left(sA+\frac{s^{3}}{3!}A^{3}+...+\frac{s^{2n+1}}{(2n+1)!}A^{2n+1}+...\right)\\
=\frac{1}{2}\left(\left(e^{s}+e^{-s}\right)E+\left(e^{s}-e^{-s}\right)\right)A.
\end{array}
\]

Hence,
\[
C(t)=\varphi(t)E+\psi(t)A,
\]
 where $\varphi(2t)=\frac{1}{2}\left(e^{t}+e^{-t}\right)$ and $\psi(2t)=\frac{1}{2}\left(e^{t}-e^{-t}\right)$.
One can check directly that
\[
\bar{A}=0,\quad\bar{B}=\frac{1}{4}\left(\begin{array}{rc}
(1+\varphi(2t)) & \quad(2+\psi(2t))\\
(2-\psi(2t)) & \quad(1-\varphi(2t))
\end{array}\right).
\]
 It is obvious that the matrix $\bar{B}$ is not constant.

\section{Preliminary study of linear systems of second-order ordinary differential equations with constant coefficients}

Let us consider a linear system of second-order ordinary differential
equations,
\begin{equation}
{\bf y}^{\prime\prime}=A{\bf y}^{\prime}+B{\bf y},\label{eq:feb27.21}
\end{equation}
 where $A$ and $B$ are matrices with constant entries. For
commuting matrices $AB=BA$ this system is reduced to the system of
the form
\begin{equation}
{\bf y}^{\prime\prime}=M{\bf y},\label{eq:feb27.22}
\end{equation}
 which is completely studied in \cite{bk:BoykoPopovychShapoval[2012],bk:WafoSoh[2010],bk:Meleshko[2011],bk:Campoamor-Stursberg[2011],bk:Campoamor-Stursberg[2012]}.
The case of noncommutative matrices is considered here.

\subsection{Equivalence transformations}

The general solution of the determining equations of the equivalence
Lie group defines the equivalence Lie group, which corresponds to
the generators
\[
X_{1}^{e}=\partial_{x},\,\,\, X_{2}^{e}=y\partial_{y}+z\partial_{z},
\]

\[
\begin{array}{rcl}
X_{3}^{e} & = & z\partial_{z}-a_{12}\partial_{a_{12}}+a_{21}\partial_{a_{21}}-b_{12}\partial_{b_{12}}+b_{21}\partial_{b_{21}},\\
X_{4}^{e} & = & y\partial_{z}-a_{12}\partial_{a_{11}}-(a_{22}-a_{11})\partial_{a_{21}}+a_{12}\partial_{a_{22}}\\
 &  & -b_{12}\partial_{b_{11}}-(b_{22}-b_{11})\partial_{b_{21}}+b_{12}\partial_{b_{22}},\\
X_{5}^{e} & = & z\partial_{y}+a_{21}\partial_{a_{11}}+(a_{22}-a_{11})\partial_{a_{12}}-a_{21}\partial_{a_{22}}\\
 &  & +b_{21}\partial_{b_{11}}+(b_{22}-b_{11})\partial_{b_{12}}-b_{21}\partial_{b_{22}},\\
X_{6}^{e} & = & x\partial_{x}-a_{11}\partial_{a_{11}}-a_{12}\partial_{a_{12}}-a_{21}\partial_{a_{21}}-a_{22}\partial_{a_{22}}\\
 &  & -2b_{11}\partial_{b_{11}}-2b_{12}\partial_{b_{12}}-2b_{21}\partial_{b_{21}}-2b_{22}\partial_{b_{22}},\\
X_{7}^{e} & = & x(y\partial_{y}+z\partial_{z})+2\partial_{a_{11}}+2\partial_{a_{22}}\\
 &  & -a_{11}\partial_{b_{11}}-a_{12}\partial_{b_{12}}-a_{21}\partial_{b_{21}}-a_{22}\partial_{b_{22}}.
\end{array}
\]

The transformations related with the generators $X_{2}^{e}$, $X_{3}^{e}$,
$X_{4}^{e}$ and $X_{5}^{e}$ correspond to the linear change of the
dependent variables $\bar{{\bf y}}=P{\bf y}$ with a constant nonsingular
matrix $P$. The generators $X_{1}^{e}$ and $X_{6}^{e}$ define shifting and scaling $x$. The transformations corresponding to the
generator $X_{7}^{e}$ define the change
\[
\begin{array}{c}
\bar{y}=ye^{\tau x},\,\,\,\bar{z}=ze^{\tau x},\,\,\,\bar{a}_{11}=a_{11}+2\tau,\,\,\,\bar{a}_{12}=a_{12},\\
\bar{a}_{21}=a_{21},\,\,\,\bar{a}_{22}=a_{22}+2\tau,\,\,\,\bar{b}_{11}=b_{11}-a_{11}\tau-\tau^{2},\\
\bar{b}_{12}=b_{11}-a_{12}\tau,\,\,\,\bar{b}_{21}=b_{21}-a_{21}\tau,\,\,\,\bar{b}_{22}=b_{22}-a_{22}\tau-\tau^{2},
\end{array}
\]
 where $\tau$ is a group parameter. In matrix form the last set of transformations
is rewritten as
\begin{equation}
\bar{{\bf y}}=e^{\tau x}{\bf y},\ \ \ \bar{A}=A+2\tau I,\,\,\,\bar{B}=B-\tau A-\tau^{2}I,\label{eq:mar10.1}
\end{equation}
where $I$ is the identical matrix.

\subsection{Simplification of the matrix $A$}

Let us apply the change $\bar{{\bf y}}=P{\bf y}$, where $P$ is a
nonsingular matrix with constant entries
\[
P=\left(\begin{array}{cc}
p_{11} & p_{12}\\
p_{21} & p_{22}
\end{array}\right).
\]
 Equations (\ref{eq:feb27.21}) become
\[
\bar{{\bf y}}^{\prime\prime}=\bar{A}\bar{{\bf y}}^{\prime}+\bar{B}\bar{{\bf y}},
\]
 where
\[
\bar{A}=PAP^{-1},\,\,\,\bar{B}=PBP^{-1}.
\]
 This means that the change $\bar{{\bf y}}=P{\bf y}$ reduces equation
(\ref{eq:feb27.21}) to the same form with the matrices $A$ and $B$
changed. Using this change, the matrix $A$ can be presented in the
Jordan form. For a real-valued $2\times2$ matrix $A$, and if the
matrix $P$ also has real-valued entries, then the Jordan matrix is
one of the following three types,
\begin{equation}
J_{1}=\left(\begin{array}{cc}
a & 0\\
0 & b
\end{array}\right),\;\; J_{2}=\left(\begin{array}{cc}
a & c\\
-c & a
\end{array}\right),\;\; J_{3}=\left(\begin{array}{cc}
a & 1\\
0 & a
\end{array}\right),\label{eq:Jordan_form-1}
\end{equation}
 where $a$, $b$ and $c>0$ are real numbers. By virtue of transformations
(\ref{eq:mar10.1}) one can assume that $a=0$. Notice also that using
the dilation of $x$, one can reduce $c$ to $1$.

In these cases one obtains the following three cases:
\[
\begin{array}{lcl}
A=J_{1|a=0} & : & AB-BA=\left(\begin{array}{cc}
0 & -b_{12}b\\
b_{21}b & 0
\end{array}\right);\\
A=J_{2|a=0,c=1} & : & AB-BA=\left(\begin{array}{cc}
(b_{12}+b_{21}) & (b_{22}-b_{11})\\
(b_{22}-b_{11}) & -(b_{12}+b_{21})
\end{array}\right);\\
A=J_{3|a=0} & : & AB-BA=\left(\begin{array}{cc}
b_{21} & (b_{22}-b_{11})\\
0 & -b_{21}
\end{array}\right);
\end{array}
\]
 where $b_{ij}$ are entries of the matrix $B$.

This then gives the following for noncommutative matrices:

(a) in the case $A=J_{1}$ one can assume that $b_{12}b\neq0$;

(b) in the case $A=J_{2}$ one has to assume that $(b_{12}+b_{21})^{2}+(b_{22}-b_{11})^{2}\neq0$;

(c) in the case $A=J_{3}$ one has to assume that $b_{21}{}^{2}+(b_{22}-b_{11})^{2}\neq0$.

\subsection{Case $A=J_{1}$}

In this case the matrix $A$ is given as follows:
\[
A=\left(\begin{array}{cc}
0 & 0\\
0 & b
\end{array}\right).
\]
 Let us denote $b=4\lambda$, where $\lambda\neq0$. The admitted generator
has the form
\[
X=C_{1}\bar{X}_{1}+C_{2}X_{2},
\]
 and the determining equations are reduced to the study of the equations
\begin{equation}
\begin{array}{c}
C_{1}(b_{11}^{3}-2b_{11}^{2}b_{22}+7b_{11}^{2}\lambda^{2}+2b_{11}b_{12}b_{21}+b_{11}b_{22}^{2}-6b_{11}b_{22}\lambda^{2}\\
-56b_{11}\lambda^{4}-2b_{12}b_{21}b_{22}-b_{22}^{2}\lambda^{2}+8b_{22}\lambda^{4}+48\lambda^{6})=0,
\end{array}\label{eq:feb_ss10}
\end{equation}
\begin{equation}
h_{1}h_{2}C_{1}=0,\label{eq:feb_ss20}
\end{equation}
\begin{equation}
(24\lambda^{4}+14b_{22}\lambda^{2}-6b_{11}\lambda^{2}-2b_{11}b_{22}+b_{12}b_{21}+2b_{22}^{2})C_{1}=0,\label{eq:feb_ss30}
\end{equation}
\begin{equation}
C_{2}b_{21}=0,\label{eq:feb_ss11}
\end{equation}
\begin{equation}
(4b_{22}+15\lambda^{2})C_{2}=0,\,\,\,(4b_{11}-\lambda^{2})C_{2}=0,\label{eq:feb_ss21}
\end{equation}
 where
\[
h_{1}=b_{11}+b_{22}+2\lambda^{2},\,\,\, h_{2}=b_{22}-b_{11}+4\lambda^{2},
\]
\[
\bar{X}_{1}=e^{-2\lambda x}\left(h_2(\partial_{x}-\lambda(y\partial_{y}-z\partial_{z}))-2\lambda b_{12}z\partial_y\right),
\]
\[
X_{2}=e^{-\lambda x}\left(2\partial_{x}-\lambda(y\partial_{y}-3z\partial_{z})\right).
\]


Let us consider equations which admit the generator $X_{2}$. In this
case $C_{2}\neq0$ and from equations (\ref{eq:feb_ss11})-(\ref{eq:feb_ss21})
one finds that
\begin{equation}
b_{11}=\lambda^{2}/4,\ \ b_{22}=-15\lambda^{2}/4,\ \ b_{21}=0,
\label{eq:mar04.1}
\end{equation}
or the matrix
\begin{equation}
\label{eq:mar15.1}
B=\frac{1}{4}
\left(\begin{array}{cc}
\lambda^2 & 4b_{12}\\
0 & -15\lambda^2
\end{array}\right).
\end{equation}
Notice that in this case $h_2=0$.
 Since in this case equations (\ref{eq:feb_ss10})-(\ref{eq:feb_ss30})
are also satisfied, then there are two admitted generators
\begin{equation}
X_{1}=e^{-2\lambda x}z\partial_{y},\,\,\, X_{2}=e^{-\lambda x}\left(2\partial_{x}-\lambda(y\partial_{y}-3z\partial_{z})\right).\label{eq:mar04.2}
\end{equation}


If the matrix $B$ is not of the form (\ref{eq:mar15.1}), then $C_2=0$, and for
the existence of a nontrivial admitted generators one needs to require that $C_1\neq 0$.
Considering equation (\ref{eq:feb_ss20}), first assume that $h_{1}=0$.
Then equations (\ref{eq:feb_ss10}) and (\ref{eq:feb_ss30}) give
\[
b_{11}=\lambda^{2},\ \ b_{22}=-3\lambda^{2},\ \ b_{21}=0.
\]
These conditions provide that $h_2=0$.
Thus the case $h_2=0$ is the general case.


Assume now that $h_{2}=0$. The general solution of equations (\ref{eq:feb_ss10})-(\ref{eq:feb_ss30})
is
\[
b_{11}=b_{22}+4\lambda^{2},\ \ b_{21}=0,
\]
or the matrix
\[
B=
\left(\begin{array}{cc}
b_{22}+4\lambda^2 & b_{12}\\
0 & b_{22}
\end{array}\right).
\]
 In this case equations (\ref{eq:feb_ss21}) are reduced to the equation
\[
C_{2}(4b_{22}+15\lambda^{2})=0.
\]

Thus, if $4b_{22}+15\lambda^{2}=0$, then one comes to the studied
case (\ref{eq:mar04.1}), where the admitted generators are (\ref{eq:mar04.2}).
If $4b_{22}+15\lambda^{2}\neq0$, then there is the only admitted
generator $X_{1}$.

In the summary the result is the following.

There are nontrivial admitted generators only if the matrix $B$ has the form:
\[
B=
\left(\begin{array}{cc}
b_{22}+4\lambda^2 & b_{12}\\
0 & b_{22}
\end{array}\right).
\]
If $4b_{22}+15\lambda^2\neq 0$, then there is the only admitted generator
\[
X_{1}=e^{-2\lambda x}z\partial_{y}.
\]
The extension of the Lie algebra, defined by the generator $X_1$, only occurs for $b_{22}=-15\lambda^2/4$.
In this case the equations admit two generators
\[
X_{1}=e^{-2\lambda x}z\partial_{y},\,\,\, X_{2}=e^{-\lambda x}\left(2\partial_{x}-\lambda(y\partial_{y}-3z\partial_{z})\right).
\]

\subsection{Cases $A=J_{2}$ and  $A=J_{3}$}

In these cases
\[
A=\left(\begin{array}{cc}
0 & 1\\
-1 & 0
\end{array}\right),\ \ A=\left(\begin{array}{cc}
0 & 1\\
0 & 0
\end{array}\right),
\]
respectively. Calculations show that there are no extensions of the trivial generators in these cases.

\section{Conclusion}

This paper completes the study of the group classification of linear
systems of two second-order ordinary differential equations with constant
coefficients. In the literature the studied case reduced to ${\bf y}^{\prime\prime}=M{\bf y}$
where $M$ is a matrix with constant entries. This study can be extended
to linear systems of second-order ordinary differential equations with more than two equations. For noncommutative matrices
we found new cases where there are additional generators in addition
to the generic ones, $\partial_{x}$ and $y\partial_{y}+z\partial_{z}$.
This implies that when considering the classification of such systems
the first step should be to check whether the matrices $A$ and $B$
commute in (\ref{eq:feb27.21}). If they commute then the case is
equivalent to the case (\ref{feb27.1}) studied in the literature
where $M$ is a constant matrix whilst if $A$ and $B$ do not commute
then the study in this paper applies and can be extended to deal with
linear systems of second-order ordinary differential equations with
more than two equations.

\section*{Acknowledgements}

SVM and SM thank the Durban University of Technology for their support
during the period of the visit and the Department of Mathematics, Statistics and Physics for their
hospitality.

\end{document}